\documentstyle[amscd,amssymb,amsthm,verbatim,amsopn,amsopn,12pt]{amsart}

\theoremstyle{plain}
\newtheorem{Thm}{Theorem}[section]

\newtheorem{Prop}[Thm]{Proposition}
\newtheorem{Cor}[Thm]{Corollary}
\theoremstyle{definition}

\newcommand{\bnum}{\begin{enumerate}}
\newcommand{\enum}{\end{enumerate}}

\newcommand{\Z}{\mathbb{Z}}

\begin{document}
\title[On non-commuting graph of a finite ring]%
	{On non-commuting graph of a finite ring}
\author[J. Dutta and D. K. Basnet]%
	{Jutirekha Dutta and Dhiren Kumar Basnet*}
\thanks{*Corresponding author}
\date{}
\maketitle
\begin{center}\small{\it Department of Mathematical Sciences, Tezpur University,}
\end{center}
\begin{center}\small{\it Napaam-784028, Sonitpur, Assam India.}
\end{center}

\begin{center}\small{\it Email: jutirekhadutta@@yahoo.com and  dbasnet@@tezu.ernet.in*}
\end{center}

\medskip

\noindent \textit{\small{\textbf{Abstract:}
The non-commuting graph $\Gamma_R$  of a finite ring $R$ with center $Z(R)$ is a simple undirected graph whose vertex set is $R \setminus Z(R)$ and two distinct vertices $a$ and $b$ are adjacent if and only if $ab \ne ba$.
In this paper,  we   show that $\Gamma_R$ is not isomorphic to certain graphs of any finite non-commutative ring $R$.  Some   connections between $\Gamma_R$ and commuting probability of $R$ are also obtained. Further, it is  shown that  the non-commuting graphs of two $\mathbb{Z}$-isoclinic rings are isomorphic if the centers of the rings have same order.  }}

\bigskip

 

\medskip

\section{Introduction}
Let $R$ be a finite ring with center $Z(R)$. The non-commuting graph of $R$, denoted by $\Gamma_R$,  is a simple undirected graph whose vertex set is $R\setminus Z(R)$ and two distinct vertices $a$ and $b$ are adjacent if and only if $ab \ne ba$.  We write $V(\Gamma_R)$ and $E(\Gamma_R)$ to denote the set of vertices    and set of edges   of $\Gamma_R$ respectively.  We also write $\deg(v)$ to denote  the degree of a vertex $v$, which is the number of edges incident on $v$. It is easy to see that
\begin{equation}\label{degree}
\deg(r) = |R| - |C_R(r)| \text{ if } r \in V(\Gamma_R)
\end{equation}
where $C_R(r) = \{x \in R : xr = rx\}$. Note that $Z(R) = \underset{r \in R}{\cap}C_R(r)$.
Many mathematicians have studied algebraic structures by means of graph theoretical properties in the last decades (see \cite{abdollahi07, abdollahi08, abdollahiakbarimaimani06, beck88, dn16, ov11} etc.).  The notion of non-commuting graph of a finite ring  was introduced by Erfanian, Khashyarmanesh and Nafar \cite{ekn15}. In Section 2, we shall show that any disconnected graph, star graph, lollipop graph or  complete bipartite graph is not isomorphic to   $\Gamma_R$  for any  finite non-commutative ring $R$. We also obtain a dominating set for  $\Gamma_R$.



The commuting probability of a finite ring $R$, denoted by $\Pr(R)$, is the probability that a randomly chosen pair of elements of $R$ commute. In 1976, MacHale \cite{machale}   introduced the concept of $\Pr(R)$. Several results on $\Pr(R)$ can be found in  \cite{BM, BMS, dB17, duttabasnetnath}. In Section 3, we obtain a formula for  $|E(\Gamma_R)|$ in terms of $\Pr(R)$ and derive several consequences. We conclude this paper showing that the non-commuting graphs of two $\Z$-isoclinic finite rings are isomorphic if the  centers of the rings have same order.

Recall that a star graph is a tree on $n$ vertices in which one vertex has degree $n - 1$ and the others have degree $1$.
A complete graph is a graph in which every pair of distinct vertices is adjacent. 
 A lollipop graph is a  graph consisting of a complete graph and a path graph, connected with a bridge. A bipartite graph is a graph whose vertex set can be partitioned into two disjoint parts in such a way that the two end vertices of every edge lie in different parts. A  complete bipartite graph is a bipartite graph such that two vertices are adjacent if and only if they lie in different parts.

\section{Some properties of $\Gamma_R$}
In this section we mainly consider the following problem: given a simple undirected graph $G$, can we find a  ring $R$ such that $\Gamma_R$ is isomorphic to $G$? In the following results, we shall show that if $G$ is a disconnected graph, star graph, lollipop graph or  complete bipartite graph then   $\Gamma_R$ is not isomorphic to $G$ for any finite ring $R$.

\begin{Prop}
Let $R$ be a finite ring. Then
\begin{enumerate}
\item  $\Gamma_R$ is connected.
\item $\Gamma_R$ is empty graph if and only if $R$ is commutative. 
\end{enumerate} 
\end{Prop}
\begin{pf}
(a) Suppose that there exists an isolated vertex $r$ in $\Gamma_R$. Then $rs = sr$ for all $s \in R$. Therefore  $r \in Z(R)$, a contradiction. Hence, the result follows.
 
  Part (b)  follows from the definition of $\Gamma_R$.
\end{pf}


\begin{Thm}\label{star}
$\Gamma_R$ is not a star or lollipop  graph for any finite non-commutative  ring $R$.
\end{Thm}
\begin{pf}
Suppose there exists a finite non-commutative ring $R$ such that $\Gamma_R$ is a star or lollipop graph. Then there exists a vertex $r$ such that $\deg(r) = 1$. This gives $[R : C_R(r)] = |R|/(|R| - 1)$, a contradiction. Hence, the result follows.
\end{pf}

In fact, the proof of the above theorem shows that  there is no vertex having degree one in $\Gamma_R$.

\begin{Thm}\label{complete relative non-commuting}
$\Gamma_R$ is not a complete bipartite graph for any finite non-commutative ring $R$. 
\end{Thm}

\begin{pf}
Let $R$ be a finite non-commutative ring such that    $\Gamma_R$ is a complete bipartite graph.
Then we have two disjoint subsets   $S_1$ and $S_2$ of $V(\Gamma_R)$ such that $|S_1| + |S_2| = |R| - |Z(R)|$. Therefore $R \cap S_1 \neq \phi$ and $R \cap S_2 \neq \phi$. Let $a \in R \cap S_1$ and $b \in R \cap S_2$. Then $ab \neq ba$. If $a + b \in R \cap S_1$ or $R \cap S_2$ then both give $ab = ba$, a contradiction. So $a + b \in Z(R)$ which gives $ab = ba$, a contradiction. Hence, the theorem follows.
\end{pf}
\begin{Thm}
$\Gamma_{R}$ is not a complete graph for any finite non  commutative ring $R$ with unity.
\end{Thm}

\begin{pf}
Let $R$ be a non-commutative ring with unity such that $\Gamma_{R}$ is complete. Then for $r \in V(\Gamma_R)$ we have 
\[
\deg(r) = |V(\Gamma_R)| - 1 = |R| - |Z(R)| - 1.
\]
By \eqref{degree}, we have $|R| - |C_R(r)| = |R| - |Z(R)| - 1$. This gives $|Z(R)| = 1$ and $|C_R(r)| = 2$, a contradiction since $|Z(R)| \geq 2$ and $|C_R(r)| \geq 3$. Hence, the result follows.
\end{pf}

We conclude this section by obtaining a dominating set for $\Gamma_R$. Recall that a dominating set of a graph $\Gamma_R$ is  a subset $D$ of $V(\Gamma_R)$ such that every vertex  in $V(\Gamma_R)\setminus D$ is adjacent to at least one member of $D$.

\begin{Prop}
Let $R$ be a finite non-commutative ring with unity. Let $A = \{a_1, a_2, \dots, a_m\}$ and $B = \{b_1, b_2, \dots, b_n\}$ be generating sets for $R$. If $A \cap Z(R) = \{a_{c + 1}, \dots, a_m\}$ and $B \cap Z(R) = \{b_{d + 1}, \dots, b_n\}$ then $D = \{a_1, a_2,\dots, a_c, b_1, b_2,\dots, b_d\}$ is a dominating set for $\Gamma_R$. 
\end{Prop}

\begin{pf}
Clearly $D \subseteq V(\Gamma_R)$. Let $r \in V(\Gamma_R)$ such that $r \notin D$. So there exists a vertex $s \in R$ such that $s = g_ib_1^{\alpha_{1i}}b_2^{\alpha_{2i}}\dots b_p^{\alpha_{pi}}$ where $g_i \in \Z$, $\alpha_{ji} \in {\mathbb{N}} \cup \{0\}$, $b_j \in B$; and $s = h_ja_1^{\alpha_{1j}}a_2^{\alpha_{2j}}\dots a_q^{\alpha_{qj}}$ where $h_j \in \Z$, $\alpha_{ij} \in {\mathbb{N}} \cup \{0\}$ and $a_i \in A$ such that $rs \neq sr$. Thus $rb_i \neq b_ir$ for some $i, 1 \leq i \leq d$ and $rs_j \neq s_jr$ for some $j, 1 \leq j \leq c$. Hence the result follows.
\end{pf}

As a corollary of the above theorem, we have the following result.

\begin{Cor}
Let $R$ be a finite non-commutative ring with unity. Let $S = \{s_1, s_2, \dots, s_n\}$ be a generating set for $R$. If $S \cap Z(R) = \{s_{m + 1}, \dots, s_n\}$ then $D = \{s_1, s_2, \dots, s_m\}$ is a dominating set for $\Gamma_R$. 
\end{Cor}

\section{Relation between $\Gamma_R$ and $\Pr(R)$}
The commuting probability of a finite ring $R$, denoted by $\Pr(R)$, is given by the following ratio
\begin{equation}\label{def_Pr}
\Pr(R) = \frac{|\{(a, b) \in R\times R : ab = ba\}|}{|R|^2}.
\end{equation} 
Note that $\Pr(R_1) = \Pr(R_2)$ if $R_1, R_2$ are two finite  non-commutative rings such that $|Z(R_1)| = |Z(R_2)|$ and $\Gamma_{R_1}, \Gamma_{R_2}$ are isomorphic. In this section, we derive the following relation between $|E(\Gamma_R)|$ and $\Pr(R)$. 
\begin{Thm}\label{relation E&R}
Let $R$ be a finite  non-commutative ring. Then the number of edges of $\Gamma_R$ is
\[
|E(\Gamma_R)| = \frac{|R|^2}{2}(1 - \Pr(R)). 
\]
\end{Thm}
\begin{pf}
Let $S = \{(a, b) \in R \times R : ab \neq ba\}$.  Then, by \eqref{def_Pr}, we have 
\begin{align*}
2|E(\Gamma_R)| = |S| &= |R|^2 - |\{(a, b) \in R \times R : ab = ba\}|\\ 
&= |R|^2 - |R|^2\Pr(R).
\end{align*}
Hence, the result follows.
\end{pf}
As a corollary to Theorem \ref{relation E&R}, we have the following lower bound for $\Pr(R)$.
\begin{Cor}
Let $R$ be a non-commutative ring. Then 
\[
\Pr(R) \geq \frac{2|Z(R)|}{|R|} + \frac{1}{|R|} -  \frac{|Z(R)|^2}{|R|^2} - \frac{|Z(R)|}{|R|^2}. 
\]  
\end{Cor}

\begin{pf}
We know that for every graph with $n$ vertices, the number of edges is at most $\frac{n(n - 1)}{2}$. Therefore
\[
|E(\Gamma_R)| \leq \frac{1}{2}(|R| - |Z(R)|)(|R| - |Z(R)| - 1).
\]
Hence, using Theorem \ref{relation E&R}, we have the required result. 
\end{pf} 
We also have the following result.
\begin{Cor}
There is no finite non-commutative ring $R$ with trivial center such that $\Pr(R) = 1 - 2/|R| + 4/|R|^2$.
\end{Cor}
\begin{pf}
Suppose there exists a finite non-commutative ring $R$ such that $|Z(R)| = 1$ and
\[
\Pr(R) = 1 - 2/|R| + 4/|R|^2.
\]
Then, by Theorem \ref{relation E&R}, we have
 \[
|E(\Gamma_R)| = |R| - |Z(R)| - 1 = |V(\Gamma_R)| -1. 
\]
This shows that there is a non-commutative ring $R$ with trivial center  such that $\Gamma_R$ is a star graph, which is a contradiction (by Theorem \ref{star}). Hence, the result follows.
\end{pf}

Now we obtain some bounds for $|E(\Gamma_R)|$ as consequences of Theorem \ref{relation E&R}.

\begin{Prop}
Let $R$ be a finite non-commutative ring and $p$ the smallest prime dividing $|R|$. Then 
\[
|E(\Gamma_R)| \leq \frac{1}{2}(|R| - |Z(R)|)(|R| - p).
\]
\end{Prop}
\begin{pf}
The result follows from  \cite[Theorem 2.1]{dB17} and Theorem \ref{relation E&R}. 
\end{pf}

\begin{Prop}
Let $R$ be a non-commutative ring. Then 
\[
|E(\Gamma_R)| \geq \frac{3|R|^2}{16}.
\]
\end{Prop}
\begin{pf}
The result follows from Theorem \ref{relation E&R} noting that $\Pr(R) \leq \frac{5}{8}$. 
\end{pf}

We conclude this section with another lower bound for $|E(\Gamma_R)|$ and a consequence of it. 
\begin{Prop}\label{lower bound for E}
Let $R$ be a finite  non-commutative ring. Then
\[
|E(\Gamma_R)| \geq \frac{|R|}{4}(|R| - |Z(R)|).
\]
\end{Prop}

\begin{pf}
By \eqref{degree}, we have
\begin{align*}
2|E(\Gamma_R)| 
&= \underset{r \in V(\Gamma_R)}{\sum}(|R| - |C_R(r)|)\\
&\geq (|R| - |Z(R)|)\left(|R| - \frac{|R|}{2}\right) = \frac{|R|}{2}(|R| - |Z(R)|).
\end{align*}
Hence we have the required result.    
\end{pf}
Using Theorem \ref{relation E&R} in Proposition \ref{lower bound for E}, we have the following upper bound for $\Pr(R)$.
\begin{Cor}\label{up2}
Let $R$ be a finite non-commutative ring. Then 
\[
\Pr(R) \leq \frac{1}{2} + \frac{1}{2}\frac{|Z(R)|}{|R|}.
\]
\end{Cor}

\section{Relation between $\Z$-isoclinism and $\Gamma_R$}

Hall \cite{pH40} introduced the notion of isoclinism between two groups and Lescot \cite{pL95} showed that the commuting probability of two isoclinic finite groups are same. Later on Buckley, MacHale and N$\acute{\rm i}$ sh$\acute{\rm e}$ \cite{BMS} introduced the concept of $\Z$-isoclinism between two rings and showed that the commuting probability of two isoclinic finite rings are same.  Recall that two rings $R_1$ and $R_2$ are said to be $\Z$-isoclinic (see \cite{BMS})  if there exist additive group isomorphisms $\phi : \frac{R_1}{Z(R_1)} \rightarrow \frac{R_2}{Z(R_2)}$   and $\psi : [R_1, R_1] \rightarrow [R_2, R_2]$ such that $\psi ([u, v]) = [u', v']$ whenever $\phi(u + Z(R_1)) = u' + Z(R_2)$ and $\phi(v + Z(R_1)) = v' + Z(R_2)$.


We have the following main result of this section.
\begin{Thm} 
Let $R_1$ and $R_2$ be two finite rings such that $|Z(R_1)| = |Z(R_2)|$. If $R_1$ and $R_2$ are $\Z$-isoclinic then $\Gamma_{R_1} \cong \Gamma_{R_2}$.  
\end{Thm}

\begin{pf}
Let $(\phi, \psi)$ be a $\Z$-isoclinism between $R_1$ and $R_2$. Then $|\frac{R_1}{Z(R_1)}| = |\frac{R_2}{Z(R_2)}|$ and $|[R_1, R_1]| = |[R_2, R_2]|$. Since $|Z(R_1)| = |Z(R_2)|$ we have  $|R_1| = |R_2|$ and $|R_1\setminus Z(R_1)| = |R_2 \setminus Z(R_2)|$. Let $\{r_1, r_2, \dots,  r_n\}$ and $\{r'_1, r'_2, \dots,   r'_n\}$  be  transversals for $\frac{R_1}{Z(R_1)}$ and  $\frac{R_2}{Z(R_2)}$ respectively.
Let $\phi$ be defined as $\phi(r_i + Z(R_1)) = r'_i + Z(R_2)$ where $r_i \in R_1$ and $r'_i \in R_2$ for all $1 \leq i \leq n$. Also, let $\theta : Z(R_1) \rightarrow Z(R_2)$ be a one-to-one correspondence. Let us define a map $\alpha : R_1 \rightarrow R_2$ such that $\alpha(r_i + z) = r'_i + \theta(z)$ for $1 \leq i \leq n$ and $z \in Z(R_1)$. Then $\alpha$ is a bijection. This shows that $\alpha$ is also a bijection from $R_1 \setminus Z(R_1)$ to $R_2 \setminus Z(R_2)$. Suppose $r, s$ are adjacent in $\Gamma_{R_1}$. Then $[r, s] \neq 0$, this gives $[r_i + z, r_j +z_1] \neq 0$ for some $z, z_1 \in Z(R_1)$, $r_i, r_j \in \{r_1, r_2, \dots,  r_n\}$ and $r = r_i + z$, $s = r_j + z_1$.

 Thus $[r'_i + \theta(z), r'_j + \theta(z_1)] \neq 0$ for some $\theta(z), \theta(z_1) \in Z(R_2)$ and $r'_i, r'_j \in \{r'_1, r'_2, \dots, r'_n\}$. Hence $[\alpha(r_i + z), \alpha(r_j + z_1)] \neq 0$, that is $\alpha(r)$ and $\alpha(s)$ are adjacent. Thus the result follows.
\end{pf}


\end{document}